\documentclass[12pt]{article}%
\usepackage{graphicx}
\usepackage{amsmath}
\usepackage{amssymb}
\usepackage{amsthm}
\usepackage{latexsym}
\usepackage{amsfonts}%
\newcommand{\comment}[1]{}
\newtheorem{theorem}{Theorem}
\newtheorem{lemma}{Lemma}[section]
\newtheorem{remark}{Remark}[section]

\newtheorem{proposition}{Proposition}[section]

\begin{document}

\bigskip\bigskip

\title{{\LARGE \textsf{Change point models and conditionally pure birth processes; an
inequality on the stochastic intensity. }}}
\date{January 30, 2004}
\author{\textbf{Emilio De Santis, Fabio Spizzichino}\\{\small \textsl{Universit\`{a} di Roma \textit{La Sapienza\/}, Dipartimento di
Matematica ``Guido Castelnuovo''}}\\{\small \textsl{Piazzale Aldo
Moro, 2 - 00185 Roma, Italia}} } \maketitle

\begin{abstract}
We analyze several aspects of a class of simple counting
processes, that can emerge in some fields of applications where
the presence of a change-point occurs. Under simple conditions we,
in particular, prove a significant inequality for the stochastic
intensity.
\end{abstract}

MSC 2000 \textit{subject Classification}: 60G55, 60J27, 60K10,
90B25.


\section{Introduction}

In this note we consider the \textit{change-point} model,
described as follows. Let a random time $U$ \ and a simple
counting process $\{N_{t} \}_{t\geq0}$ be defined on a same
probability space. Let $T_{1}\leq T_{2} \leq\dots$ be the arrival
times of $\{N_{t}\}_{t\geq0}$ and denote by $h_{t}$ an
\textquotedblleft history", observed in the time-interval $[0,t]$,
such as
\[
h_{t}\equiv\{T_{1}=t_{1},\dots,T_{k}=t_{k},T_{k+1}>t\}
\]
with $0\equiv t_{0}<t_{1}<t_{2}<\dots<t_{k}<t$.

We think of $U$ as a change-point for $\{N_{t}\}_{t\geq0}$; more
precisely we assume that the latter admits an intensity, described
by
\begin{equation}
\lim_{\Delta\rightarrow0^{+}}\frac{1}{\Delta}P(T_{k+1}<t+\Delta\,|\,h_{t}
;\,U\leq t)=\lambda_{1}(k) \label{int1}
\end{equation}
\begin{equation}
\lim_{\Delta\rightarrow0^{+}}\frac{1}{\Delta}P(T_{k+1}<t+\Delta\,|\,h_{t}
;\,U>t)=\lambda_{0}(k), \label{int2}
\end{equation}
where $\lambda_{0}(0),\lambda_{0}(1),\dots$ and
$\lambda_{1}(0),\lambda _{1}(1),\dots$ are two given sequences of
positive constants; i.e., given the observation of the history
$h_{t}$, $\lambda_{1}(\cdot)$ would be the intensity, conditional
on the knowledge that time $t$ is after the change-point $U$ and
$\lambda_{0}(\cdot)$ would be the intensity, conditional on the
knowledge that time $t$ is before $U$. We think of the case when
the random variable $U$ is not observable (obviously $U$ takes
values in the interval $(0, \infty)$); its distribution function
will be denoted by $G$.

More formally, let $X_{t}\equiv1_{\{U<t\}}$ and denote by
$\Im\equiv \{{\mathcal{F}}_{t}^{(N)}\}_{t\geq0}$ and
$\aleph\equiv\{{\mathcal{F}}_{t}^{(N,X)}\}_{t\geq0}$ the
filtrations respectively generated by the processes
$\{N_{t}\}_{t\geq0}$ and $\{N_{t},X_{t}\}_{t\geq0}$; we are then
assuming that the stochastic intensity of $\{N_{t}\}_{t\geq0}$
with respect to $\aleph$ is
\[
\label{sjs}\mu_{t}^{\aleph}=\lambda_{1}(N_{t})X_{t}+\lambda_{0}(N_{t}
)(1-X_{t}).
\]
The intensity of $\{N_{t}\}_{t\geq0}$ w.r.t. the ``internal"
filtration $\Im$ is then specified by the position:
\[
\label{inten}\mu_{t}^{\Im}=\lambda_{1}(N_{t})P(U\leq t|{\mathcal{F}}_{t}%
^{(N)})+\lambda_{0}(N_{t})P(U>t|{\mathcal{F}}_{t}^{(N)});
\]
we shall also use the following notation:
\begin{equation}
\mu_{t}(h_{t})=\lambda_{1}(k)P(U\leq t\,|\,h_{t})+\lambda_{0}%
(k)P(U>t\,|\,h_{t}). \label{int3}%
\end{equation}
A counting process that is a Pure Birth process conditionally on
the change point, as described so far, will be denoted by the
symbol CPB$(G,\lambda _{0}(\cdot),\lambda_{1}(\cdot))$.

Assume now
\begin{equation}
\label{assu}\lambda_{1} (k ) \geq\lambda_{0} ( k ) \hbox{ for } k
= 0,1, \dots
\end{equation}
and compare two different observed histories
\begin{equation}
\label{ist1}h^{\prime}_{t} \equiv\{ T_{1} = t^{\prime}_{1}, \dots,
T_{k} = t^{\prime}_{k} , T_{k+1} >t\}
\end{equation}
\begin{equation}
\label{ist2}h^{\prime\prime}_{t} \equiv\{ T_{1} =
t^{\prime\prime}_{1}, \dots, T_{k} = t^{\prime\prime}_{k} ,
T_{k+1} >t\}
\end{equation}
both containing $k$ arrivals in the same time-interval $[0,t]$.

We write $h_{t}^{\prime\prime}\unrhd h_{t}^{\prime}$ if
\begin{equation}
t_{i}^{\prime\prime}\geq t_{i}^{\prime}\hbox{ for }i=1,2,\dots
,k\label{napoli0}%
\end{equation}
or, equivalently, for any $s\in\lbrack0,t)$,
\begin{equation}
\sum_{j=1}^{\infty}1_{\{s\leq t_{j}^{\prime\prime}\leq t\}}=\sum_{j=1}%
^{k}1_{\{s\leq t_{j}^{\prime\prime}\leq
t\}}\geq\sum_{j=1}^{k}1_{\{s\leq t_{j}^{\prime}\leq
t\}}=\sum_{j=1}^{\infty}1_{\{s\leq t_{j}^{\prime}\leq
t\}}\label{napoli}%
\end{equation}
i.e. $h_{t}^{\prime\prime}\unrhd h_{t}^{\prime}$ when the number
of recent arrivals in the history $h_{t}^{\prime}$ is not larger
than the number of recent arrivals in $h_{t}^{\prime\prime}$.

In the present note we analyze some aspects of CPB counting
processes and, in particular, we prove the following result.

\begin{theorem}
\label{p1} If $h_{t}^{\prime\prime}\unrhd h_{t}^{\prime}$ then $\mu_{t}%
(h_{t}^{\prime\prime})\geq\mu_{t}(h_{t}^{\prime})$.
\end{theorem}

Our interest for this result is illustrated in the following
remark.

\begin{remark}
\label{r2} When $U$ is exponentially distributed, the computation
of $\mu _{t}(h_{t})$ can be in principle carried out explicitly;
in fact, in such a case, one can compute the normalizing constant
that is needed to obtain the conditional probabilities in Eq.
(\ref{int3}) and this in turn allows Theorem~\ref{p1} to admit a
direct proof. The case with $\lambda
_{1}(0)=\lambda_{1}(1)=\cdots$ and
$\lambda_{0}(0)=\lambda_{0}(1)=\cdots$ is dealt with in
\cite{HJ02}; more lengthy expressions may be involved in our case
where $\lambda_{0}$ and $\lambda_{1}$ may depend on the number of
past arrivals. The explicit computation of the normalizing
constant in Eq.~(\ref{int3}) is however not possible when $U$ is
not exponentially distributed.

Obviously $\{N_{t}\}_{t\geq0}$ is not a pure birth-process (i.e.
it is not Markov): when we ``uncondition" with respect to the
random variable ${\mathbf{1}}_{\{U\leq t\}}$ in Eq. (\ref{int3}),
we obtain an intensity $\mu_{t}(h_{t})$ which depends on the
arrival times $t_{1},t_{2},\dots,t_{k}$ and not only on $k$. It is
then natural to wonder whether it is possible to
establish some \textit{a priori } inequalities on the pair $(\mu_{t}%
(h_{t}^{\prime}),\mu_{t}(h_{t}^{\prime\prime}))$.
\end{remark}

The paper will be organized as follows. In Sect. 2 we will
consider a random change of time-scale that will reveal to be
useful in the proof of Theorem~\ref{p1}; we will in particular
show that the class of CPB processes is closed under this type of
transformation.

Theorem~\ref{p1} will be proved in Section~\ref{s2}. On this
purpose, we prove an analogous result for a corresponding
discrete-time model, afterwards the desired result will be
obtained by means of a suitable passage to the limit. We notice
that the discrete time result can be however of autonomous
interest.

Section \ref{s4} will be devoted to a brief discussion and to some
final remarks on Theorem~\ref{p1} and on the class of CPB counting
processes. Models in this class emerge in a natural way in several
fields; in particular we shall mention two cases of interest, in
the frame of reliability and experimental sciences, respectively.

For several aspects of the well-known change-point problem and a
comprehensive bibliography, we address the reader \textit{e.g.} to
\cite{Ar, C, HJ02, MS93} and references therein; we refer to
Bremaud \cite{Br1, Br2} for general aspects about counting
processes. For properties of monotonicity and of stochastic
orderings for counting processes, see \cite{KS96} and \cite{MS02}.

\bigskip

\section{A random time-scale transformation}

Besides the process $\{N_{t} \}_{t\geq0}$, we shall introduce in
this section a new counting process
$\{{\tilde{N}}_{t}\}_{t\geq0}$; Lemma~\ref{lemsec} and
Proposition~\ref{le1} to be obtained below will turn out to be
useful for our purposes in the next section. Lemma~\ref{lemsec} in
particular shows that the conditional probability of the event
$\{U>t\}$, given an observed history $h_{t}$ for
$\{{N}_{t}\}_{t\geq0}$, does coincide with an analogous
conditional probability for $\{{\tilde{N}}_{t}\}_{t\geq0}$.

Such a new process, which also admits intensities, is obtained
from the original one by means of a random change of time-scale,
as follows.

Let $(\Omega,{\mathcal{F}},P)$ be the probability space on which
the random variables $U,T_{1},T_{2},\dots$ are defined and let
$\gamma_{0},\gamma _{1},\dots$ be a sequence of positive constants
with
\[
\label{bound1}0<\inf_{i}\gamma_{i}\leq\sup_{i}\gamma_{i}< \infty.
\]
Let $g:\Omega\times\lbrack0,\infty)\rightarrow\lbrack0,\infty)$ be
the strictly increasing random function of time defined as follows
\begin{equation}
\label{tyhg}g(\omega, t)=\sum_{k=0}^{N_{t}(\omega)-1}\gamma_{k}(T_{k+1}%
(\omega)- T_{k}(\omega))+\gamma_{N_{t}(\omega)}(t-T_{N_{t}(\omega)
}(\omega)) .
\end{equation}
>From now on the symbol $\omega$ will be dropped; using a more compact
notation we write
\[
g(t) =
\int_{0}^{t}\sum_{k=0}^{\infty}\gamma_{k}{1}_{[T_{k},T_{k+1})}
(s)ds,
\]
or, by setting
\[
\label{nuob}\gamma(s)=\sum_{k=0}^{\infty}\gamma_{k}1_{[T_{k},T_{k+1})}(s),
\]
\[
\label{nuob2}g(t)=\int_{0}^{t}\gamma(s)ds.
\]

Define now, on $(\Omega, {\mathcal{F}}, P)$, the random variables
\begin{equation}
\label{ese1}{\tilde U} = g( U )=\sum_{k=0}^{N_{U} - 1} \gamma_{k}
(T_{k+1} -T_{k}) +\gamma_{N_{U}} (U -T_{N_{U}}) =\int_{0}^{U}
\gamma(s ) ds
\end{equation}
\begin{equation}
\label{ese2}{\tilde T}_{l} = g( T_{l})= \sum_{k=0}^{l - 1}
\gamma_{k} (T_{k+1} -T_{k}) = \int_{0}^{T_{l}} \gamma(s ) ds
\hbox{ for } l = 1,2 , \dots.
\end{equation}
and consider the new counting process $\{{\tilde N}_{t} \}_{t
\geq0 }$ whose arrival times are ${\tilde T}_{1}$, ${\tilde
T}_{2}, \dots$; thus we have
\begin{equation}
\label{nuovoP}{\tilde N}_{g (t ) } =N_{t} .
\end{equation}

Let $A_{1},A_{2},...$ denote the \textit{interarrival times} of $\{N_{t}%
\}_{t\geq0}$:
\[
N_{t}=\sup\{n|\sum_{k=1}^{n}A_{k}\leq t\}.
\]

Notice that the transformation yielding
$\{\widetilde{N}_{t}\}_{t\geq0}$ can also be described by writing
\[
\widetilde{N}_{t}=\sup\{n|\sum_{k=1}^{n}\gamma_{k-1}A_{k}\leq t\},
\]
i.e. $\{ \widetilde{N}_{t} \}_{t \geq0 }$ is such that its
interarrival times $\widetilde{A}_{1}, \,\, \widetilde{A}_{2},
\dots$, satisfy
\begin{equation}
\label{ultim}\widetilde{A}_{k} = \gamma_{k-1 } A_{k}.
\end{equation}

We denote ${\tilde X}_{t} \equiv1_{\{ {\tilde U} < t \}}$, and
consider the filtrations ${ \tilde\Im} \equiv\{ {
{\mathcal{F}}}^{({\tilde N})}_{t}\}_{t \geq0 }$, ${ \tilde\aleph}
\equiv\{ { {\mathcal{F}}}^{({\tilde N} , {\tilde X})}_{t}\}_{t
\geq0 }$. From now on, for typographic convenience, we shall often
use the symbols $N(t)$ and ${\tilde N} (t)$ in place of $N_{t} $
and ${\tilde N }_{t}$, respectively.

As we shall see, the interest in the transformation defined by
(\ref{tyhg}), is motivated by the following Lemma.

\begin{lemma}
\label{lemsec} Under the positions (\ref{tyhg}), (\ref{ese1}), and
(\ref{ese2}) one has

\begin{itemize}
\item[a)] ${{\mathcal{F}}}^{({\tilde N} , {\tilde X})}_{g (t)} =
{{\mathcal{F}}}^{({ N} , { X})}_{t}$ and ${
{\mathcal{F}}}^{({\tilde N} )}_{g (t)} = {{\mathcal{F}}}^{({ N}
)}_{t}$;

\item[b)] $
P(U > t | {{\mathcal{F}}}^{({ N} )}_{t}) = P({\tilde U} > g(t) | {
{\mathcal{F}}}^{({\tilde N} )}_{g (t)})$.
\end{itemize}
\end{lemma}

\textbf{Proof.} First we notice the following: since the
transformation $g$ defined by (\ref{tyhg}) is continuous and
increasing in $t$, we have
\[
\{ \omega\in\Omega: N_{t} = k, T_{1} \in[ t_{1} , t_{1}+
\Delta_{1}), \dots, T_{k} \in[ t_{k} , t_{k}+ \Delta_{k} ) ,
T_{k+1} >t , U >s\} =
\]
\begin{equation}
\label{fregene}\big \{ \omega\in\Omega: {\tilde N}_{g(t)} = k ,
{\tilde T}_{1} \in{\big [}g( t_{1}) , g(t_{1}+ \Delta_{1}) {\big
)}, \dots, {\tilde T}_{k} \in{\big [} g(t_{k}) ,g( t_{k}+
\Delta_{k}) {\big )}, {\tilde T}_{k+1} >g(t) , {\tilde U} >g(s)
\big \} ,
\end{equation}
and
\[
\big \{ \omega\in\Omega: N_{t} = k, T_{1} \in\big [ t_{1} , t_{1}+
\Delta _{1}\big ), \dots, T_{k} \in\big [ t_{k} , t_{k}+
\Delta_{k} \big ) , T_{k+1}
>t \big \} =
\]
\begin{equation}
\label{civita}= \big \{ \omega\in\Omega: {\tilde N}_{g(t)}=k,
{\tilde T}_{1} \in\big [ g( t_{1}) , g(t_{1}+ \Delta_{1}) \big ),
\dots, {\tilde T}_{k} \in\big [ g(t_{k}) ,g( t_{k}+ \Delta_{k})
\big), {\tilde T}_{k+1} >g(t) \big \} .
\end{equation}
a) ${{\mathcal{F}}}^{({ N} , { X})}_{t}$ is actually generated by
the subsets of the type
\[
\{ \omega\in\Omega: N_{t} = k, T_{1} \in[ t_{1} , t_{1}+
\Delta_{1}), \dots, T_{k} \in[ t_{k} , t_{k}+ \Delta_{k}) ,
T_{k+1} >t , U >s\}
\]
with $k = 0,1 , \dots$, $s \leq t$, $0 \leq t_{1} \leq\cdots\leq
t_{k} \leq t$,
\newline and ${ {\mathcal{F}}}^{({\tilde N} , {\tilde X})}_{g (t)}$ is generated by
the subsets of the type
\[
\big \{ \omega\in\Omega: { {\tilde N}_{g(t)}=k, {\tilde T}_{1}
\in{\big [}g( t_{1}) , g(t_{1}+ \Delta_{1}) {\big )}, \dots,
{\tilde T}_{k} \in{\big [} g(t_{k}) ,g( t_{k}+ \Delta_{k}) {\big
)}, {\tilde T}_{k+1} >g(t) , \tilde U} >g(s) \big \} .
\]
Then the identity ${ {\mathcal{F}}}^{({\tilde N} , {\tilde X})}_{g
(t)} = {{\mathcal{F}}}^{({ N} , { X})}_{t}$ follows from
(\ref{fregene}). Similarly ${ {\mathcal{F}}}^{({\tilde N} )}_{g
(t)} = {{\mathcal{F}}}^{({ N} )}_{t}$ follows from (\ref{civita}).

b) The assertion immediately follows from a) by noticing that $\{U
> t\} = \{{\tilde U} > g(t)\} $. \fbox{}

\begin{remark}
\label{remmm} The intensity $\mu_{t} ( h_{t}) $ can be expressed
in terms of the process $\{\widetilde{N}_{t}\}_{t\geq0} $; in
fact, by b) of Lemma~\ref{lemsec} and by recalling the notation
(\ref{int3}), we can write
\begin{equation}
\label{perteo}\mu_{t} ( h_{t})= \lambda_{0} (k) + (\lambda_{1} (k
)-\lambda_{0} (k)) P({\tilde U} \leq g(t) \, | \,{\tilde T}_{1} =
g(t_{1}) , \dots, {\tilde T}_{k}= g(t_{k}), {\tilde T}_{k+1} >
g(t) ) .
\end{equation}

\end{remark}

The class of CBP processes is closed under the transformation
defined by (\ref{tyhg}); more precisely we have the following
result that can also be inspired by equation (\ref{ultim}).

\begin{proposition}
\label{le1} The process $\{{\tilde{N}}_{t} \}_{t\geq0}$ is CPB$({\tilde{G}%
},{\tilde{\lambda}}_{1},{\tilde{\lambda}}_{0})$ where
${\tilde{G}}(s)\equiv P({\tilde{U}}<s)$ and
\begin{equation}
{\tilde{\lambda}}_{i}(k)=\frac{\lambda_{i}(k)}{\gamma_{k}}. \label{formula1}%
\end{equation}

\end{proposition}

Actually Proposition \ref{le1} could be proved by using a general,
well known, result about simple counting processes; the latter
shows how a simple counting process, admitting intensity, can be
obtained from a standard Poisson process via a random change of
time scale (see e.g.\cite{Br2, Ku}). For the reader's convenience,
we prefer however to give here a direct proof which uses the
specific notation of this paper.

The following remark will be used in such a proof.

\begin{remark}
\label{remm} The events $\{ {\tilde N}(g({ t}+ \Delta) ) > {\tilde
N}({ g(t)}) \} $ and $\{ {\tilde N}(g(t)+ \gamma_{{\tilde N}
(g(t))}\Delta)> {\tilde N } (g(t)) \}$ are equal. In order to show
such identity we notice:
\begin{equation}
\label{rmf}\{ N( t+ \Delta) > N(t) \} = \{ {\tilde N}(g({ t}+
\Delta) )
>{\tilde N}({ g(t)}) \}
\end{equation}
furthermore, if there is no arrival in the interval $( t , t +
\Delta]$ for the original process $N$, we can write
\begin{equation}
\label{rmf2}g( t + \Delta) = g(t ) + \gamma_{{\tilde N}
(g(t))}\Delta= g(t ) + \gamma_{N (t)}\Delta,
\end{equation}
whence ${\tilde N} (g( t + \Delta) ) = {\tilde N} ( g(t ) +
\gamma_{N (t)}\Delta)$, i.e. we can conclude that, if $N(t +
\Delta) = N (t )$ then ${\tilde N} (g(t ) + \gamma_{N (t)}\Delta)=
{\tilde N} (g(t)) $ as well.

Let us suppose, on the other hand, that there are one or more
arrivals for the process $N$ in the interval $(t, t + \Delta]$ and
denote by $T_{a}$ the instant of the earliest among these such
arrivals. Then there is an arrival for ${\tilde N } $ at the
instant $g(T_{a})$, that is within the interval $(g(t ) , g(T_{a})
] = ( g(t ) , g(t) + \gamma_{{\tilde N} (g(t))} (T_{a} - t ) ] $.
This means that there is at least one arrival in $( g(t) , g(t)
\gamma_{{\tilde N} (g(t))} \Delta]$; in fact $( g(t) , g(t)
\gamma_{{\tilde N} (g(t))} \Delta] \supseteq( g(t ) , g(t) +
\gamma_{{\tilde N} (g(t))} (T_{a} - t ) ]$, since
$\Delta\geq(T_{a} -t )$.
\end{remark}

\textbf{Proof of Proposition \ref{le1}.} We know that
\begin{equation}
\lim_{\Delta\rightarrow0^{+}}\frac{P(N({t}+\Delta)-N({t})>0|{\mathcal{F}}%
_{t}^{(N,X)})}{\Delta}=\lambda_{X(t)}(N_{t}), \label{fi1}%
\end{equation}
and, taking into account Lemma \ref{lemsec} and Eq.
(\ref{nuovoP}),
\begin{equation}
\lim_{\Delta\rightarrow0^{+}}\frac{P(N(t+\Delta)-N(t)>0|{\mathcal{F}}%
_{t}^{(N,X)})}{\Delta}=\lim_{\Delta\rightarrow0^{+}}\frac{P({\tilde{N}}%
(g({t}+\Delta))-{\tilde{N}}({g(t)})>0|{{\mathcal{F}}}_{g(t)}^{({\tilde{N}%
},{\tilde{X}})})}{\Delta}. \label{ds1}%
\end{equation}

On the other hand in view of the previous Remark~\ref{remm} we
have
\[
\lim_{\Delta\to0^{+} }\frac{P( {\tilde N}(g({ t}+ \Delta) ) -
{\tilde N}({ g(t)}) >0 | {{\mathcal{F}}}^{({\tilde N},{\tilde
X})}_{g(t)})}{ \Delta} =\lim_{\Delta\to0^{+} }\frac{P( {\tilde
N}(g(t)+ \gamma_{{\tilde N}(g(t))}
\Delta) - {\tilde N } (g(t)) >0 | { {\mathcal{F}}}^{({\tilde N},{\tilde X}%
)}_{g(t)} )}{ \Delta}
\]
\[
=\lim_{\Delta\to0^{+} }\frac{P( {\tilde N}(g(t)+ \gamma_{{ \tilde N}%
(g(t))}\Delta)- {\tilde N } (g(t)) >0 | { {\mathcal{F}}}^{({\tilde
N},{\tilde X})}_{g(t)} )}{ \gamma_{{\tilde N}(g(t))}\Delta}
\frac{\gamma_{{\tilde N}(g(t))}\Delta}{\Delta}
\]

\[
\label{fi2c}=\gamma_{{\tilde{N}}(g(t))}\lim_{\epsilon\rightarrow0^{+}}%
\frac{P({\tilde{N}}(g(t)+\epsilon)-{\tilde{N}}(g(t))>0|{{\mathcal{F}}}%
_{g(t)}^{({\tilde{N}},{\tilde{X}})})}{\epsilon}.
\]
Then, by (\ref{fi1}), we can write
\[
\label{fi5}\lim_{\epsilon\rightarrow0^{+}}\frac{P({\tilde{N}}(g(t)+\epsilon
)-{\tilde{N}}(g(t))>0|{{\mathcal{F}}}_{g(t)}^{({\tilde{N}},{\tilde{X}})}%
)}{\epsilon}=\frac{\lambda_{X(t)}(N_{t})}{\gamma_{{\tilde{N}}(g(t))}}%
=\frac{\lambda_{X(t)}(N_{t})}{\gamma_{{N}(t)}}%
\]
i.e.
\[
\label{fi2d}{\tilde{\lambda}}_{\tilde{X}(g(t))}({\tilde{N}}_{g(t)}%
)=\lim_{\epsilon\rightarrow0^{+}}\frac{P({\tilde{N}}(g(t)+\epsilon)-{\tilde
{N}}(g(t))>0|{{\mathcal{F}}}_{g(t)}^{({\tilde{N}},{\tilde{X}})})}{\epsilon
}=\frac{\lambda_{X(t)}(N_{t})}{\gamma_{{N}(t)}}.
\]
Thus for $i=0,1$ and for every $k=0,1,\dots$ we obtain
\[
\label{pqq}{\tilde{\lambda}}_{i}(k)=\frac{\lambda_{i}(k)}{\gamma_{k}}.
\]
This completes the proof. \fbox{}

\section{\label{s2} Discrete approximations and proof of Theorem \ref{p1}}

We will start this section by considering a discrete approximation
of the continuous-time model; this will allow us to prove a
discrete-time version of Theorem~\ref{p1} under an additional
condition (see (\ref{ser}) below).

Afterwards, by performing a natural limit, we will obtain the
desired result for the continuous-time model. In order to
eliminate the condition (\ref{ser}) we shall resort to the
counting process $\{{\tilde{N}}_{t}\}_{t\geq0}$ and to the related
results obtained in the previous section.

Consider a discrete-time model defined as follows. Let $\bar{U}$
be an ${\mathbb{N}}$-valued random time and set, for $m = 1,2,
\dots$
\begin{equation}
{\bar{X}}_{n}\equiv1_{\{{\bar{U}}<n\}},\hbox{ }\nu(m)\equiv P({\bar{U}%
}=m|{\bar{U}}>m-1), \label{capalbio}%
\end{equation}
so that
\begin{equation}
P( \bar{U} =1) = \nu(1) , \,\,\,\,\,\, P({\bar{U}}=m)=\nu(m)\prod_{l=1}%
^{m-1}(1-\nu(l)), \label{ansedonia}%
\end{equation}
and we assume
\begin{equation}
\label{orbetello}\nu(m)\in(0,1) .
\end{equation}

Let ${\bar T}_{1} < {\bar T}_{2} < \dots$ be an increasing
sequence of ${\mathbb{N}}$-valued random times and set
\begin{equation}
\label{lki}{\bar N}_{n} \equiv\sup\{ k | {\bar T}_{k} \leq n\}.
\end{equation}
We assume that two sequences of positives constants
$\{{\bar\lambda}_{0}(k ) \}_{k=0,1 , \dots}$ and
$\{{\bar\lambda}_{1}(k ) \}_{k=0,1 , \dots}$ exist such that
\begin{equation}
\label{lki2}%
\begin{array}
[c]{c}%
P( {\bar T_{k+1} =n+1 | {\bar T}_{1} = n_{1} , \dots, {\bar T}_{k}
= n_{k}
,{\bar T}_{k+1} >n, {\bar U} >n }) = {\bar\lambda}_{0}(k )\\
P( {\bar T_{k+1} =n+1 | {\bar T}_{1} = n_{1} , \dots, {\bar T}_{k}
= n_{k} , {\bar T}_{k+1} >n,{\bar U} \leq n }) = {\bar\lambda}_{1}
(k )
\end{array}
\end{equation}
and that, for any $k \in{\mathbb{N}}$
\begin{equation}
\label{plo}{\bar\lambda}_{1} ( k ) >{\bar\lambda}_{0} ( k ) .
\end{equation}
Furthermore we assume here
\begin{equation}
\label{ser}\frac{(1- {\bar\lambda}_{1} (k-1)) (1-
{\bar\lambda}_{0} (k))} {(1- {\bar\lambda}_{0} (k-1)) (1-
{\bar\lambda}_{1} (k))} \geq1 \hbox { for } k = 1 , \dots.
\end{equation}
For an history ${\bar h}_{n} \equiv\{ {\bar T}_{1} = n_{1}, \dots,
{\bar T}_{k} = n_{k}, {\bar T}_{k+1} > n \} $, (where $0\leq n_{1}
< \dots< n_{k} \leq n $) we set
\[
\label{grosseto}{\bar\mu}_{n} ({\bar h}_{n}) \equiv P( {\bar
T}_{k+1} = n+1 | {\bar h}_{n} )= {\bar\lambda}_{1} (k ) P({\bar U
} \leq n \, | \, {\bar h}_{n} ) + {\bar\lambda}_{0} (k ) P({\bar U
} > n \, | \, {\bar h}_{n} ) .
\]

In Proposition~\ref{propf} below we will use the following
observation

\begin{remark}
\label{az} Let $A$, $B$, $C$, $D$, $\alpha$, $\beta$, $\gamma$ be
positive constants and define
\begin{equation}
\theta= \frac{C}{A+B+C+ D},
\end{equation}
\begin{equation}
\theta^{\prime}= \frac{C \gamma}{A\alpha+B \gamma+C\gamma+ D
\delta}.
\end{equation}
If $\alpha/ \gamma\geq1$ and $\delta/ \gamma\geq1$ then
$\theta\geq \theta^{\prime}$ .

\end{remark}

\begin{proposition}
\label{propf} Under the conditions (\ref{orbetello}), (\ref{plo})
and (\ref{ser}), we have
\begin{equation}
{\bar{\mu}}_{n}({\bar{h}^{\prime\prime}}_{n})\geq{\bar{\mu}}_{n}({\bar
{h}^{\prime}}_{n}) \label{livorno}%
\end{equation}
for any pair ${\bar{h}^{\prime}}_{n}$,
${\bar{h}^{\prime\prime}}_{n}$ such that
${\bar{h}}_{n}^{\prime\prime}\unrhd{\bar{h}}_{n}^{\prime}$.
\end{proposition}

\textbf{Proof}. First we notice that the inequality
(\ref{livorno}) is equivalent to
\begin{equation}
{P}({\bar{U}}>n\,\,|\,\,{\bar{h}}_{n}^{\prime\prime}\,)\leq{P}({\bar{U}%
}>n\,\,|\,{\bar{h}}_{n}^{\prime}).\label{did}%
\end{equation}
Now we denote, for $1\leq n_{1}<n_{2}<\cdots<n_{k}\leq n$ and
$j=1,2,\dots$
\[
\label{cassia}g_{j}(n_{1},n_{2},\dots,n_{k};n)=P({\bar{U}}=j,{\bar{T}}%
_{1}=n_{1},\dots,{\bar{T}}_{k}=n_{k},{\bar{T}}_{k+1}>n),
\]
so that
\begin{equation}
{P}({\bar{U}}>n\,\,|\,\,{\bar{h}}_{n})=\frac{\sum_{j=n+1}^{\infty}g_{j}%
(n_{1},n_{2},\dots,n_{k};n)}{\sum_{j=1}^{\infty}g_{j}(n_{1},n_{2},\dots
,n_{k};n)}.\label{disc1}%
\end{equation}
In view of (\ref{ansedonia}) and (\ref{lki2}),
\[
g_{j}(n_{1},n_{2},\dots,n_{k};n)=
\]%
\begin{equation}
\left[  \nu(j)\prod_{r=1}^{j-1}{(1-\nu(r))}\right]
\prod_{i=0}^{k-1}\left\{
{\bar{\lambda}}_{{\mathbf{1}}(n_{i+1}>j)}(i)\prod_{r=n_{i}+1}^{n_{i+1}%
-1}[1-{\bar{\lambda}}_{{\mathbf{1}}(r>j)}(i)]\right\}  \prod_{r=n_{k}+1}%
^{n}[1-{\bar{\lambda}}_{{\mathbf{1}}(r>j)}(i)]\label{disc2b}%
\end{equation}
where we have set $n_{0}=0$.

Now, on the space of possible ``discrete'' histories, let us
define the operators $\Phi_{i}$ as follows. For an history ${\bar
h}_{n} \equiv\{ {\bar T}_{1} = n_{1}, \dots, {\bar T}_{i} =n_{i},
\dots, {\bar T}_{k} = n_{k}, {\bar T}_{k+1} > n \} $, let
\[
\label{pisa1}\Phi_{i} ({\bar h}_{n}) \equiv\{ {\bar T}_{1} =
n_{1}, \dots, {\bar T}_{i} =n_{i}+1, \dots, {\bar T}_{k} = n_{k},
{\bar T}_{k+1} > n \}
\]
for $i $ such that $i= 1, \dots, k-1 $ and $n_{i+1} > n_{i} +1$,
\[
\label{pisa2}\Phi_{k} ({\bar h}_{n}) \equiv\{ {\bar T}_{1} =
n_{1}, \dots, {\bar T}_{k} = n_{k}+1, {\bar T}_{k+1} > n \}
\,\,\,\,\, \hbox{ if } n_{k} < n,
\]
and
\[
\label{pisa3}\Phi_{i} ({\bar h}_{n}) \equiv{\bar h}_{n}%
\]
otherwise.

It is easy to check that any history ${\bar h^{\prime\prime}}_{n}
$ such that ${\bar h}^{\prime\prime}_{n} \unrhd{\bar
h}^{\prime}_{n} $ can be obtained from ${\bar h}^{\prime}_{n}$, by
applying the operators $\Phi_{i}$ a finite number of times.

Then we can reduce ourselves to show the validity of the
inequality
\begin{equation}
\label{veron}{ P} ( {\bar U} > n \,\, | \,\, {\bar h}_{n} ) \geq{
P} ( {\bar U} > n \,\, | \,\, \Phi_{l} ({\bar h}_{n} ) ) \hbox{
for } l \in\{ 1,2 \dots, k \} .
\end{equation}

With $l $ as in (\ref{veron}) we now let
\[
\label{denum1}A( n_{1} , \dots, n_{k} ;n )= \sum_{j=1}^{n_{l}
-1}g_{j} (n_{1} , \dots, n_{k}; n),
\]
\[
\label{denum2}B( n_{1} , \dots, n_{k} ;n )=
\sum_{j=n_{l}+1}^{n}g_{j} (n_{1} , \dots, n_{k}; n) ,
\]
\[
\label{num}C( n_{1} , \dots, n_{k} ;n )=
\sum_{j=n+1}^{\infty}g_{j} (n_{1} , \dots, n_{k};n ).
\]
Then we can rewrite formula (\ref{disc1}) as
\[
\label{tuscolana}{ P} ( {\bar U} > n \,\, | \,\, {\bar h}_{n} ) =
\frac{C( n_{1} , \dots, n_{k} ;n)}{ A( n_{1} ,\dots, n_{k};n )+
g_{n_{l}} (n_{1} , \dots, n_{k};n) +B( n_{1} , \dots, n_{k};n )+C(
n_{1} , \dots, n_{k} ;n) } .
\]
We now switch to obtaining the expression of ${ P} ( {\bar U} > n
\,\, | \,\, \Phi_{l} ({\bar h}_{n} ) ) $ in terms of $A( n_{1},
\dots, n_{k} ; n )$, $B( n_{1}, \dots, n_{k} ; n )$, $C( n_{1},
\dots, n_{k} ; n )$ and $g_{n_{l} }( n_{1}, \dots, n_{k} ; n )$.

Let us denote
\[
\label{viterbo}%
\begin{array}
[c]{c}%
{\hat{A}}(n_{1},\dots,n_{k};n)=\sum_{j=1}^{n_{l}-1}g_{j}(n_{1},\dots
,n_{l}+1,\dots,n_{k};n),\\
{\hat{B}}(n_{1},\dots,n_{k};n)=\sum_{k=n_{l}+1}^{n}g_{k}(n_{1},\dots
,n_{l}+1,\dots,n_{k};n),\\
{\hat{C}}(n_{1},\dots,n_{k};n)=\sum_{k=n+1}^{\infty}g_{k}(n_{1},\dots
,n_{l}+1,\dots,n_{k};n),\\
{\hat{g}}={g}_{n_{l}}(n_{1},\dots,n_{l}+1,\dots,n_{k};n).
\end{array}
\]
The following identities hold:
\begin{equation}
\label{ma}{\hat{A}}(n_{1},\dots,n_{k};n)=\alpha\cdot
A(n_{1},\dots,n_{k};n),
\end{equation}
\begin{equation}
\label{mb}{\hat{B}}(n_{1},\dots,n_{k};n)=\gamma\cdot
B(n_{1},\dots,n_{k};n),
\end{equation}
\begin{equation}
\label{mc}{\hat{C}}(n_{1},\dots,n_{k};n)=\gamma\cdot
C(n_{1},\dots,n_{k};n),
\end{equation}
\begin{equation}
\label{md}{\hat{g}}=g_{n_{l}}(n_{1},\dots,n_{l}+1,\dots,n_{k};n)=\delta\cdot
g_{n_{l}}(n_{1},\dots,n_{l},\dots,n_{k};n),
\end{equation}
with
\[
\label{viterbo2}\alpha=\frac{1-{\bar{\lambda}}_{1}(l-1)}{1-{\bar{\lambda}}%
_{1}(l)},\hbox{ }\delta=\frac{[1-{\bar{\lambda}}_{0}(l-1)]{\bar{\lambda}}%
_{1}(l-1)}{{\bar{\lambda}}_{0}(l-1)[1-{\bar{\lambda}}_{1}(l)]},\hbox{
}\gamma
=\frac{1-{\bar{\lambda}}_{0}(l-1)}{1-{\bar{\lambda}}_{0}(l)}.
\]

\medskip

In order to check the validity of the identity (\ref{ma}),\ we can
just notice that, for $1\leq j\leq n_{l}-1$, it is
\[
g_{j}(n_{1},\dots,n_{l}+1,\dots,n_{k};n)=\frac{1-{\bar{\lambda}}_{1}%
(l-1)}{1-{\bar{\lambda}}_{1}(l)}g_{j}(n_{1},\dots,n_{l},\dots,n_{k};n),
\]
in view of formula (\ref{disc2b}); then
\[
{\hat{A}}(n_{1},\dots,n_{k};n)=\sum_{j=1}^{n_{l}-1}g_{j}(n_{1},\dots
,n_{l}+1,\dots,n_{k};n)=
\]

\[
=\frac{1-{\bar{\lambda}}_{1}(l-1)}{1-{\bar{\lambda}}_{1}(l)}\sum_{j=1}%
^{n_{l}-1}g_{j}(n_{1},\dots,n_{l},\dots,n_{k};n)=\alpha\cdot A(n_{1}%
,\dots,n_{k};n).
\]

The validity of the identities (\ref{mb})-(\ref{md}) can be
obtained in an analogous way with $j>n_{l}$ or $j=n_{l}$,
respectively.

By the definitions of ${\hat{A}}$, ${\hat{B}}$, ${\hat{C}}$ and
${\hat{g}}$ and by taking into account the identities
(\ref{ma})-(\ref{md}), we can now write
\[
\label{casilina}{P}({\bar{U}}>n\,\,|\,\,\Phi_{l}({\bar{h}}_{n}))=\frac
{{\hat{C}}}{{\hat{A}}+{\hat{g}}+{\hat{B}}+{\hat{C}}}=\frac{\gamma
C}{\alpha A+\delta g+\gamma B+\gamma
C}=\frac{C}{\frac{\alpha}{\gamma}A+\frac{\delta }{\gamma}g+B+C}.
\]
 It is immediately seen (by using Remark \ref{az}) that, in view
of the assumptions (\ref{plo})-(\ref{ser}), it is
\[
{P}({\bar{U}}>n\,\,|\,\,\Phi_{l}({\bar{h}}_{n}))\leq{P}({\bar{U}%
}>n\,\,|\,\,{\bar{h}}_{n})
\]
and this proves the assertion. \fbox{}

We are now in a position to prove Theorem~\ref{p1}.

\textbf{Proof.} Let us assume for the moment that, besides the
condition (\ref{assu}), also the following condition holds:
\begin{equation}
\label{catania}\lambda_{1} ( k ) - \lambda_{0} ( k) < \lambda_{1}
( k+1 ) - \lambda_{0} ( k+1) \hbox{ for } k=0,1 \dots
\end{equation}

Consider a sequence of discrete-time models as follows: for $m =
1,2, , \dots$ let ${\bar U^{(m)}} \equiv\frac{[ m U ]}{ m} $ and
${\bar T}^{(m)}_{l} $ ($l= 1,2 \dots$) be discrete random
variables taking values on the set $\{ 0, \frac{1}{m},\frac{2}{m},
\dots\}$ and such that
\begin{equation}
\label{viareggio}{\bar\lambda}^{(m)}_{0}(k) \equiv\frac{{ \lambda}_{0}(k )}%
{m}, \,\,\,\,\,\, {\bar\lambda}^{(m)}_{1}(k) \equiv\frac{{
\lambda}_{1}(k )}{m},
\end{equation}
where we set
\begin{equation}
\label{viarnord}%
\begin{array}
[c]{c}%
{\bar\lambda}^{(m)}_{0}(k) \equiv P( {{\bar T}^{(m)}_{k+1}
=\frac{n+1}{m} | {\bar T}^{(m)}_{1} = \frac{n_{1}}{m} , \dots,
{\bar T}^{(m)}_{k} = \frac
{n_{k}}{m} ,{\bar T}^{(m)}_{k+1} >\frac{n}{m}, {\bar U} >\frac{n}{m} }) ,\\
{\bar\lambda}^{(m)}_{1}(k) \equiv P( {{\bar T}^{(m)}_{k+1}
=\frac{n+1}{m} | {\bar T}^{(m)}_{1} = \frac{n_{1}}{m} , \dots,
{\bar T}^{(m)}_{k} = \frac {n_{k}}{m} ,{\bar T}^{(m)}_{k+1}
>\frac{n}{m}, {\bar U} \leq\frac{n}{m} }) .
\end{array}
\end{equation}

It can be checked that, for $h_{t} =\{ T_{1}= t_{1} , \dots, T_{k}
=t_{k} , T_{k+1} > t\}$
\begin{equation}
\label{lido}P(U >t | h_{t} ) = \lim_{m \to\infty} P \left(  {\bar
U}^{(m)} > \frac{[tm]}{m} {\Big |} {\bar T}_{1} =\frac{[t_{1}
m]}{m}, \dots, {\bar T}_{k+1} > \frac{[tm]}{m} \right)  .
\end{equation}
We do not report all the details; we limit ourselves to mention
that, in order to obtain the identity (\ref{lido}), one has first
to take into account
\begin{equation}
\label{marina}P( U > t | h_{t} ) = \lim_{\Delta\to0 } \frac{P( U>t
, T_{1} \in[ t_{1} , t_{1} + \Delta) , \dots, T_{k} \in[ t_{k} ,
t_{k} + \Delta) , T_{k+1} > t )}{ P( T_{1} \in[ t_{1} , t_{1} +
\Delta) , \dots, T_{k} \in[ t_{k} , t_{k} + \Delta) , T_{k+1} > t
) } .
\end{equation}
The r.h.s. of (\ref{marina}) can be shown to be equal to
\[
\label{forte}\lim_{m \to\infty} \frac{\sum_{n =0}^{\infty} P(
T_{1} \in[ t_{1} , t_{1} + \Delta_{m}) , \dots, T_{k} \in[ t_{k} ,
t_{k} + \Delta_{m} | A^{(m)}_{n} (t))) P(A^{(m)}_{n}(t))}{\sum_{n
=0}^{\infty} P( T_{1} \in[ t_{1} , t_{1} + \Delta_{m}) , \dots,
T_{k} \in[ t_{k} , t_{k} + \Delta_{m} |
A^{(m)}_{n}(0) )) P(A^{(m)}_{n}(0)) }%
\]
where we have denoted, for $s \geq0 $, $A^{(m)}_{n} (s)=\{ {\bar
U}^{(m)} \in[ s + n\Delta_{m} , s + (n+1) \Delta_{m} ) \}$ and
$\Delta_{m} $ is an infinitesimal sequence.

Finally the identity (\ref{lido}) can be obtained by a
Poisson-type approximation by taking into account the position
(\ref{viareggio}); for a general discussion about Poisson
approximations and for results similar to the one needed here see
e.g. \cite{BHJ}.

Assuming (\ref{assu}) and (\ref{catania}) we obtain the conditions
(\ref{plo})-(\ref{ser}) for the intensities
${\bar\lambda}^{(m)}_{0}(k) $ and ${\bar\lambda}^{(m)}_{1}(k) $,
in fact the condition (\ref{plo}) is trivially verified for all
the integer $m $ and the condition (\ref{ser}) is easily obtained
from (\ref{catania}) by a Taylor expansion, for large $m $.
Consider now $h^{\prime}_{t}$ and $h^{\prime\prime}_{t}$ as in
(\ref{ist1}), (\ref{ist2}) and, for $m $ large enough, the
corresponding histories in discrete time defined by
\[
\label{massa}%
\begin{array}
[c]{c}%
{\bar h}^{\prime}_{t} =\{ {\bar T}_{1} =\frac{[t^{\prime}_{1}
m]}{m}, \dots,
{\bar T}_{k+1}>\frac{[tm]}{m} \}\\
{\bar h}^{\prime\prime}_{t} =\{ {\bar T}_{1}
=\frac{[t^{\prime\prime}_{1} m]}{m}, \dots, {\bar
T}_{k+1}>\frac{[tm]}{m} \} .
\end{array}
\]
If $h^{\prime\prime}_{t} \unrhd h^{\prime}_{t} $ then, it is
\begin{equation}
\label{carrara}\frac{[t^{\prime}_{i} m]}{m}
\leq\frac{[t^{\prime\prime}_{i} m]}{m} \hbox{ for } i=1, \dots, k
..
\end{equation}
In view of condition (\ref{carrara}) we have the inequality
(\ref{did}) (see the proof of Proposition~\ref{propf}); by using
(\ref{lido}), we can then obtain, for the continuous-time limit
process $\{N_{t} \}_{t \geq0 }$ the inequality
\begin{equation}
\label{didcont}{P}(U>t \,\,|\,\,{{h}}_{t}^{\prime\prime}\,)\leq{P}%
(U>t\,\,|\,h_{t}^{\prime}).
\end{equation}
whence $\mu_{t}(h_{t}^{\prime\prime})\geq\mu_{t}(h_{t}^{\prime})$
immediately follows under the assumption (\ref{catania}).

We now show however that such an assumption is by no means
restrictive. Suppose in fact that we deal with a
CPB$(G,\lambda_{1}, \lambda_{0})$ with $\lambda_{1}, \lambda_{0}$
not satisfying (\ref{catania}) and apply the transformation $g(t)$
in (\ref{tyhg}) with the specific choice
\[
\label{pescara}\gamma_{k} = c_{k} \frac{\lambda_{1} (k ) -
\lambda_{0} (k ) }{ \lambda_{1} (0 ) - \lambda_{0} (0 ) } .
\]
Here $c_{k}$ is an arbitrary decreasing sequence such that $c_{0}
=1$ and $\lim_{k \to\infty} c_{k} >0 $. By Lemma \ref{le1} we thus
obtain a new CPB$({\tilde G},{\tilde\lambda}_{1} ,
{\tilde\lambda}_{0})$ process where
\[
\label{tirreno}{\tilde\lambda}_{i} (k) = \frac{\lambda_{i} (k)}{
\gamma_{k}} = \lambda_{i} (k)\frac{1}{c_{k}}\frac{\lambda_{1} (0 )
- \lambda_{0} (0 ) }{
\lambda_{1} (k ) - \lambda_{0} (k ) }%
\]
so that the condition (\ref{catania}) is satisfied.
Lemma~\ref{lemsec} then shows that we were entitled to prove
Theorem~\ref{p1} for the CPB$({\tilde G},{\tilde\lambda}_{1} ,
{\tilde\lambda}_{0})$. \fbox{}

\medskip\medskip\medskip

Obviously the conditional probability $P(U< t | h_{t}) $ in
(\ref{int3}) actually depends on the parameters $\lambda_{i} ( k )
$ ($i = 0,1 $ and $k = 0,1, \dots$); then it might be convenient
to use the notation $P_{\{\lambda _{i} (\cdot)\}}(U< t | h_{t}) $.

Notice that the inequality in (\ref{assu}) has been taken in the
strict sense; however, by using the continuity property of
$P_{\{\lambda_{i} (k)\}}(U< t | h_{t}) $ with respect to the set
of parameters $\lambda_{i} (\cdot)$, it can be easily checked that
in Theorem \ref{p1} the inequality can be taken in the broad
sense.

\section{\label{s4} Discussion and concluding remarks.}

The notion of CPB processes, as it has been described in the
Introduction, is a very natural model, that can emerge in several
fields of application; it can be used to formalize a number of
possible situations, that, apart from the use of different
languages, turn out to be substantially isomorphic one to the
other. Here we give just two possible instances, taken from
different fields of application.

{ \textbf{Example} \textit{(A reliability application)}}. A
typical problem in reliability modelling is the description of
stochastic dependence among lifetimes of components that are to
operate simultaneously in a same environment; two simple models of
dependence in this respect, are quite common in the reliability
literature: the standard \textit{change-point }model and the
\textit{load-sharing} model.

The standard change-point model can be described as follows: $n$
components $C_{1}, \dots, C_{n}$, that we assume to be identical
for simplicity's sake, start operate simultaneously and go on
working, each $C_{i}$ until its own failure time $W_{i}$ and with
no physical interaction with the others. However $C_{1},...,C_{n}$
are imbedded in a same environment and it is the case that the
environmental condition will suddenly change its state at a random
time $U$ (the change point); this creates a form of stochastic
dependence among the
failure times $W_{1},...,W_{n}$: conditionally on $\{U=u\}$, $W_{1}%
,...,W_{n}\,$ are independent with a same failure rate coinciding
with a given function $\rho_{0}\left(  t\right)  $ for $t<u$ and
coinciding with a different failure rate function $\rho_{1}\left(
t\right)  $ for $t\geq u.$

The load-sharing model emerges instead when $C_{1},...,C_{n}$
share a same load or share a benefit from a same favorable
external condition: this makes that, between two subsequent
failure times $W_{(i)}$ and $W_{(i+1)}$, the components that
survived $W_{(i)}$ act independently, with a failure rate function
dependent on the overall number $(n-i)$ of surviving components
and, possibly, on the calendar time.

This situation is described by the fact that the counting process
$\{ N_{t}\}_{t \geq0 }$ with
\[
N_{t}=\sum_{i=1}^{n}1_{\{T_{i}\leq t\}}%
\]
is Markov (possibly non-homogeneous), i.e.\@ it is a pure death
process; for more details on this aspect see e.g. \cite{Ar2} and
\cite{SSS02}; for some examples and a wider list of references on
the load-sharing model, see also \cite{Ro84}, \cite{ShS91},
\cite{Sp}.

The CPB models considered in the present paper arise as a natural
superposition of standard change-point and load-sharing models, as
described so far. In fact, conditionally on the change point $U$,
the failure-times $W_{1}, \dots, W_{n}$ are not independent, but
rather they obey a common load-sharing model and the counting
process $\{ N_{t}\}_{t\geq0 }$ is a CPB process. This is of
interest in that one may often have to handle sets of components
that share the same load (or the same stress) and the latter can
suddenly increase its level at an unpredictable instant. The
question may arise in those cases whether, under a same number of
observed failures within a time-instant $t$, we have to be more
pessimistic with early failure times or with very recent failure
times.

Theorem 1 gives a response to this question under the condition
that in any case the hazard of surviving components becomes more
severe after the change-point.

{ \textbf{Example} \textit{(An application in Physics)}}. Under a
different language, the same superimposition of a change- point
model and a load-sharing model, can be of interest in the field of
experimental sciences. One can think for instance of $n$ spins
$C_{1},...,C_{n}$ embedded into a uniform magnetic field;
initially all the spins are in the state $-1$ and each of them
flips to its ground state $+1$ in a random time $W_{i}$. We assume
that, at any time-instant, the transition rate, beside being an
increasing function of the intensity of the magnetic field, is
influenced also by the number of already flipped spins.

Furthermore we think of the cases where the intensity of the
magnetic field, at time $0$, has the value $B_{0}$ and, at a
random time-instant $U$ flips to value $B_{1}$, with
$B_{1}>B_{0}$. The CPB model applies when the underlying magnetic
field is not directly observable.

\bigskip

Of course the examples above concern the case of counting
processes with a finite number of arrivals in the interval
$[0,\infty)$; examples of interest also can be found for the case
of infinite arrivals.

We now conclude the paper with a remark about the pair of
histories to be compared.

In Theorem \ref{p1} we compared two histories observed on the same
time-interval $\left[  0,t\right]  $ and containing the same
number of arrivals $k$. Consider now two different histories
$h_{t}^{\prime}$ and $h_{t}^{\prime\prime}$ on the same
time-interval $\left[  0,t\right]  $ where $h_{t}^{\prime\prime}$
is obtained from $h_{t}^{\prime}$ by simply "adding" some
arrivals. Under assumption (\ref{assu}), one may guess that the
inequality $\mu_{t}
(h_{t}^{\prime\prime})\geq\mu_{t}(h_{t}^{\prime})$ holds.

We notice on the contrary that this is not true, as the following
simple example shows. Let $\{N_{t}\}_{t\geq0}$ be a CPB$\left(
G,\lambda_{0}\left(
\cdot\right)  ,\lambda_{1}\left(  \cdot\right)  \right)  $ process with%
\[
\overline{G}\left(  t\right)  =\exp\{-t\},
\]
\[
\lambda_{0}\left(  0\right)  =\lambda_{0}\left(  1\right)
=1;\lambda _{1}\left( 0\right)  =2,\lambda_{1}\left(  1\right)
=M\gg2.
\]
and simply consider the two histories
\[
h_{t}^{\prime}\equiv\{\text{no arrival in }\left[  0,t\right]  \};
\,\,\,\,\,\,\, h_{t}^{\prime\prime}\equiv\{T_{1}=t_{t},T_{2}>t\}.
\]
It is easy to check that $\mu_{t} (h_{t}^{\prime\prime})< \mu_{t}%
(h_{t}^{\prime})$.

Some more assumptions on the conditional birth rates are then
needed in order to get the inequality
$\mu_{t}(h_{t}^{\prime\prime})\geq\mu_{t}(h_{t}^{\prime })$.

\bigskip

Some considerations analogous to those above can be made
concerning the comparison between two histories that contain the
same number of arrivals but are observed over two different
time-intervals: consider e.g. the two histories
$h_{t^{\prime}}^{\prime}$ and
$h_{t^{\prime\prime}}^{\prime\prime}$ be given by
\[
h_{t^{\prime}}^{\prime} = \{ \hbox{ no arrivals in } [0,
t^{\prime}]\} ; \,\,\, h_{t^{\prime\prime}}^{\prime\prime} = \{
\hbox{ no arrivals in } [0, t^{\prime\prime}] \}
\]
with $t^{\prime}< t^{\prime\prime}$. It is clear that with
appropriate choice of $G$ and of the rates $\lambda_{i} (0)$
($i=0,1$) we can have
\[
P(U > t^{\prime}| h_{t^{\prime}}^{\prime}) \leq P(U >
t^{\prime\prime}| h_{t^{\prime\prime}}^{\prime\prime})
\]
or
\[
P(U > t^{\prime}| h_{t^{\prime}}^{\prime}) > P(U >
t^{\prime\prime}| h_{t^{\prime\prime}}^{\prime\prime}).
\]
We can then conclude that, by only assuming the condition
(\ref{assu}), an inequality as in Theorem~\ref{p1} cannot be
obtained by comparing two histories if they are not observed on
the same time-interval and do not contain the same number of
arrivals.

\end{document}